\documentclass{gtmon_a}
\pdfoutput=1

\usepackage[all]{xy}
\usepackage{stmaryrd}

\proceedingstitle{Proceedings of the Nishida Fest (Kinosaki 2003)}
\conferencestart{28 July 2003}
\conferenceend{8 August 2003}
\conferencename{International Conference in Homotopy Theory}
\conferencelocation{Kinosaki, Japan}

\editor{Matthew Ando}
\givenname{Matthew}
\surname{Ando}

\editor{Norihiko Minami}
\givenname{Norihiko}
\surname{Minami}

\editor{Jack Morava}
\givenname{Jack}
\surname{Morava}

\editor{W Stephen Wilson}
\givenname{W Stephen}
\surname{Wilson}

\title{Toward a fundamental groupoid\\for the stable homotopy category}

\author{Jack Morava}
\givenname{Jack}
\surname{Morava}
\address{Department of Mathematics\\
Johns Hopkins University\\\newline
Baltimore MD 21218\\
USA}
\email{jack@math.jhu.edu}
\urladdr{}

\volumenumber{10}
\issuenumber{}
\publicationyear{2007}
\papernumber{17}
\startpage{293}
\endpage{317}

\doi{}
\MR{}
\Zbl{}

\arxivreference{math/0509001}


\received{28 May 2005}
\revised{9 July 2006}
\accepted{}
\published{18 April 2007}
\publishedonline{18 April 2007}
\proposed{}
\seconded{}
\corresponding{}
\version{}

\let\xysavmatrix\xymatrix
\def\xymatrix{\disablesubscriptcorrection\xysavmatrix}
\AtBeginDocument{\let\bar\wbar\let\tilde\wtilde\let\hat\what}

\newtheorem{conj}{Conjecture}

\theoremstyle{remark}
\newtheorem{subsect}{\hspace{-4pt}}[section]
\newtheorem{subsubsect}{\hspace{-4pt}}[subsect]

\makeautorefname{subsect}{\hspace{-4pt}}
\makeautorefname{subsubsect}{\hspace{-4pt}}

\newcommand{\cA}{\mathcal{A}}
\newcommand{\aut}{\mathrm{Aut}}
\newcommand{\e}{\mathbf{e}}
\newcommand{\F}{\mathbb{F}}
\newcommand{\bF}{\wbar{\mathbb F}}
\newcommand{\BF}{\wbar{F}}
\newcommand{\fF}{\mathfrak{F}}
\newcommand{\ff}{\mathfrak{f}}
\renewcommand{\fg}{\mathfrak{g}}
\newcommand{\G}{\mathbb{G}}
\newcommand{\cG}{\mathcal{G}}
\newcommand{\Gal}{\mathrm{Gal}}
\newcommand{\cH}{\mathcal{H}}
\newcommand{\Hom}{\mathrm{Hom}}
\newcommand{\bI}{\mathbf{I}}
\newcommand{\bk}{\wbar{k}}
\newcommand{\hK}{\hat{K}}
\newcommand{\MZN}{\mathrm{MZN}}
\newcommand{\oh}{\mathbf{o}}
\newcommand{\fp}{\mathfrak{p}}

\newcommand{\bpi}{\pmb{\pi}}
\newcommand{\QSymm}{\mathrm{QSymm}}
\newcommand{\bQ}{\wwbar{\mathbb{Q}}}
\makeop{spec}
\newcommand{\cT}{\mathcal{T}}
\newcommand{\Vect}{\mathrm{Vect}}
\newcommand{\bV}{\mathbb{V}}
\newcommand{\W}{\mathbb{W}}

\newcommand{\w}{\mathrm{W}}
\newcommand{\hZ}{\what{\mathbb{Z}}}
\newcommand{\ab}{\mathrm{ab}}
\newcommand{\hot}{\mathrm{hot}}
\newcommand{\mot}{\mathrm{mot}}

\newcommand{\odd}{\mathrm{odd}}
\newcommand{\nr}{\mathrm{nr}}
\newcommand{\bnr}{\mathbf{nr}}
\newcommand{\tr}{\mathrm{tr}}

\begin{document}

\begin{abstract}
This very speculative sketch suggests that a theory of fundamental
groupoids for tensor triangulated categories could be used to describe
the ring of integers as the singular fiber in a family of ring-spectra
parametrized by a structure space for the stable homotopy category,
and that Bousfield localization might be part of a theory of `nearby'
cycles for stacks or orbifolds.
\end{abstract}

\maketitle

\section*{Introduction}
\label{sec0}

One of the motivations for this paper comes from John Rognes' Galois
theory for structured ring
spectra. His account \cite{37} ends with some very interesting remarks about
analogies between
classical primes in algebraic number fields and the non-Euclidean primes
of the stable homotopy
category, and I try here to develop a language in which these analogies
can be restated as the
assertion that Waldhausen's \emph{unfolding}
\[
\spec \Z \to \spec S
\]
of the integers in the category of brave new rings (or $E_\infty$
ring-spectra, or commutative
$S$--algebras) leads to the existence of commutative diagrams of the form
\[\xymatrix{
{\spec \oh_{\bQ_p}}\ar[r]\ar[d]^{\Gal(\bQ_p/\Q_p)}&{\spec L^{MU}_{K(n)}MU}
\ar[d]^{D^\times}\\
{\spec \Z_p} \ar[r]&{\spec L_{K(n)}S \;.}
}\]
The vertical arrow on the left is the Galois cover defined by the ring
of integers in an algebraic
closure of the $p$--adic rationals, but to make sense of the right-hand
side would require, among
other things, a good theory of structure objects (analogous to the prime
ideal spectra of
commutative algebra) for some general class of tensor triangulated
categories. In this direction I
have mostly hopes and analogies, summarized in \fullref{sec3}.

The main result of the first section below, however, is that local
classfield theory implies the
existence of an interesting system of group homomorphisms
\[
\rho \co \w(\bQ_p/\Q_p) \to D^\times
\]
which can plausibly be interpreted as the maps induced on the fundamental
groups of these hypothetical structure objects. [The group on the left
is Weil's technical variant of the Galois group; its topology (see
Lichtenbaum~\cite{30}) is slightly subtler than the usual one.]

This example is local: it depends on a choice of the prime $p$. But Rognes
\cite[Section~12.2.1]{37} has
global results as well; in particular, he identifies the ringspectrum
$S[BU]$ as the Hopf algebra of
functions on an analog $\Gal_\hot$ of a Galois group for $MU$, regarded
as an (inseparable)
algebraic closure of $S$. This seems to be in striking agreement with
work of Connes and
Marcolli~\cite{8}, Cartier~\cite{5}, and Deligne and Goncharov~\cite{12}
on a remarkable `cosmic' generalization of Galois theory, involving
a certain motivic Galois group $\Gal_\mot$: there is a very natural
morphism
\[
S[BU] \to H\Z \otimes \QSymm_*
\]
of Hopf algebra objects, $\QSymm_*$ being a certain graded Hopf algebra
of quasisymmetric functions (see Hazelwinkel~\cite{19}), defined
by the inclusion of the symmetric in the quasisymmetric functions,
and in \fullref{sec2} I suggest that a quotient of this map defines
a representation
\[
\Gal_\mot \to \Gal_\hot
\]
(at least, over $\Q$) which conjecturally  plays the role of the
homomorphism induced on fundamental groups by a `geometric realization'
construction \cite{33}, sending the derived category of mixed Tate
motives to some category of complex-oriented spectra.

The paper ends with a very impressionistic (fauvist?) discussion of a
possible theory of fundamental groupoids for (sufficiently small) tensor
triangulated categories, which might be flexible enough to encompass
both these examples.

I am indebted to many mathematicians for conversations about this material
over the years, most recently A~Ne'eman, T~Torii, and W~Dwyer. It is
a particular pleasure to dedicate this paper to Goro Nishida, in
recognition of his broad vision of the importance of group actions,
at many levels, in topology.

\section{Some local Galois representations}
\label{sec1}

This section summarizes some classical local number theory. The first
two subsections define a
system $\w (\bQ_p/\Q_p) \to D^\times$ of representations of certain
Galois-like groups in the
units of suitable $p$--adic division algebras; then \fullref{sec1.3} sketches a
conjecture about the structure
of the absolute inertia group $\Gal(\bQ_p/\Q_p^\nr)$, which says, roughly,
that these systems are
nicely compatible.

\addtocounter{subsect}{-1}
\begin{subsect}
A \emph{local} field is a commutative field, with a nontrivial
topology in which it is
locally compact. The reals or complexes are examples, but I will be
concerned here mostly with
totally disconnected cases, in particular the fields of characteristic
zero obtained as
non-Archimedean completions of algebraic number fields. These are finite
extensions $L$ of
$\Q_p$, for some Euclidean prime $p$; the topology defines a natural
equivalence class of
valuations, with the elements algebraic over $\Z_p$ as (local) valuation
ring.

The Galois group $\Gal(\bQ/\Q)$ of an algebraic closure of the rationals
acts on the set of
prime ideals in the ring of algebraic integers in $\Q$, with orbits
corresponding to the classical
primes; the corresponding isotropy groups can be identified with the
Galois groups
$\Gal(\bQ_p/\Q_p)$ of the algebraic closures of the $p$--adic rationals.
These isotropy groups
preserve the valuation rings, and hence act on their residue fields,
defining (split) exact sequences
\[
1 \to \bI(\bQ_p/\Q_p) \to \Gal(\bQ_p/\Q_p) \to \Gal(\bF_p/\F_p) \cong
\hZ \to 0
\]
of profinite topological groups. The cokernel is the closure of a dense
subgroup $\Z$ generated
by the Frobenius automorphism $\sigma\co x \mapsto x^p$ of $\bF_p$, and
Weil observed that in
some contexts it is more natural to work with the pullback extension
\[
1 \to \bI(\bQ_p/\Q_p) \to \w(\bQ_p/\Q_p) \to \Z \to 0 \;,
\]
the so-called Weil group, which is now only locally compact. [These
groups are defined much more generally by Tate~\cite{41} but we won't
need that here.]

\addtocounter{subsect}{1}
\begin{subsubsect}
\label{sec1.1.1}
The work of Lubin and Tate defines very interesting
representations of these groups
as automorphisms of certain one-dimensional formal groups. Honda's
logarithm
\[
\log_q(T) = \sum_{i \geq 0} p^{-i} T^{q^i}
\]
(with $q = p^n$) defines a formal group law
\[
F_q(X,Y) = \log_q^{-1}(\log_q(X) + \log_q(Y)) \in \Z_p\llbracket X,Y
\rrbracket
\]
(it's not obvious that its coefficients are integral!) whose endomorphism
ring contains, besides
the elements
 \[
T \mapsto [a](T) = \log_q^{-1}(a \log_q(T))
\]
defined by multiplication in the formal group by a $p$--adic integer $a$,
the endomorphism
\[
T \mapsto [\omega](T) = \log_q^{-1}(\omega \log_q(T)) = \omega T
\]
defined by multiplication by a $(q-1)$st root of unity $\omega$. In fact
the full ring of
endomorphisms of $\F_q$ can be identified with the Witt ring $W(\F_q)$,
which can also be
described as the algebraic integers in the unramified extension field
$\Q_q$ obtained from
$\Q_p$ by adjoining $\omega$.
\end{subsubsect}

\begin{subsubsect}
By reducing the coefficients of $F_q$ modulo $p$ we obtain
a formal group law
$\BF_q$ over $\F_p$, which admits
\[
T \mapsto F(T) = T^p
\]
as a further endomorphism. Honda's group law is a Lubin--Tate group (see
Serre~\cite{39}) for the field $\Q_q$, and it can be shown that
\[
[p](T) \equiv T^q \; {\rm mod} \; p \;;
\]
in other words, $F^n = p$ in the ring of endomorphisms of $\BF_q$. From
this it is not hard to
see that
\[
{\rm End}_{\bF_p}(\BF_q) = W(\F_q) \langle F \rangle/(F^n - p)
\]
(the pointed brackets indicating a ring of noncommutative indeterminates,
subject to the relation
\[
a^\sigma F = F a
\]
for $a \in W(\F_q)$, with $\sigma \in \Gal(\F_q/\F_p) \cong \Z/n\Z$
being the (other!)
Frobenius endomorphism). This is the ring $\oh_D$ of integers in the
division algebra
\[
D = \Q_q \langle F \rangle/(F^n - p)
\]
with center $\Q_p$; its group $\oh_D^\times$ of strict units is the full
group of automorphisms
of $\BF_q$, and it is convenient to think of
\[
1 \to \oh_D^\times \to D^\times \to \Z \to 0
\]
as a semidirect product, with quotient generated by $F$. Thus conjugation
by $F$ acts on
$W(\F_q)^\times \subset \oh_D^\times$ as $\sigma \in \Gal(\F_q/\F_p)$:
the Galois action is
encoded in the division algebra structure.
\end{subsubsect}

\addtocounter{subsect}{1}
\setcounter{subsubsect}{0}
\begin{subsubsect}
\label{sec1.2.1}
For each $n \geq 1$, a deep theorem of Weil and Shafarevich
defines a
homomorphism
\[
ws_q \co \w(\Q_q^\ab/\Q_p) \cong W(\F_q)^\times \ltimes \Z \to \oh_D^\times
\ltimes \Z =
D^\times
\]
of locally compact topological groups. It extends the local version of
Artin's reciprocity law,
which defines an isomorphism
\[
L^\times \cong \w(L^\ab/L)
\]
for any totally disconnected local field $L$ (with $L^\ab$ its maximal
abelian extension). When
$L = \Q_q$ is unramified, the sequence
\[
1 \to \w(\Q_q^\ab/\Q_q) \to \w(\Q_q^\ab/\Q_p) \to \w(\Q_q/\Q_p) \to 0
\]
is just the product
\[
0 \to W(\F_q)^\times \times \Z \to W(\F_q)^\times \ltimes \Z \to \Z/n\Z
\to 0
\]
of the elementary exact sequence
\[
0 \to \Z \to \Z \to \Z/n\Z \to 0
\]
with a copy of the units in $W(\F_q)$: in the semidirect product extension
above, the generator
$1 \in \Z$ acts on $W(\F_q)$ by $\sigma$, so its $n$th power acts
trivially.
\end{subsubsect}

\begin{subsubsect}
More generally, any Galois extension $L$ of degree $n$
over $\Q_p$ can be
embedded in a division algebra of rank $n$ with center $\Q_p$, as a
maximal commutative
subfield. When the division algebra $D$ has invariant $1/n$ in the Brauer
group $\Q/\Z$ of
$\Q_p$ (as it does in our case: this invariant equals the class, modulo
$\Z$, of the $p$--order of
an element (eg $F$) generating the maximal ideal of $\oh_D$), the
theorem \cite[Appendix]{45} of
Weil and Shafarevich defines an isomorphism of the normalizer of
$L^\times$ in $D^\times$ with
the Weil group of $L^\ab$ over $\Q_p$.

[The maximal abelian extension $L^\ab$ is obtained by adjoining the
$p$--torsion points of the
Lubin--Tate group of $L$ to the maximal unramified extension $L^\nr$;
the resulting field
acquires an action of the group $\oh_L^\times$ of units of $L$ (by
`complex multiplication'
on the Lubin--Tate group) together with an action of the automorphisms
$\hZ$ of the algebraic
closure of the residue field.

The point is that the Lubin--Tate group of $L$ is \emph{natural} in the
etale topology: to be
precise, any two Lubin--Tate groups for $L$ become isomorphic over the
completion $L^\bnr$
of the maximal unramified extension $L^\nr = L \otimes_{W(k)} W(\bk)$
of $L$ (see Serre~\cite[Section~3.7]{39}). Since an
automorphism of $L$ over $\Q_p$ takes one Lubin--Tate group to another,
the resulting group
of automorphisms of `the' Lubin--Tate group of $L$ (as a completed Hopf
algebra over $\Z_p$)
is an extension of $\Gal(L/\Q_p)$ by $\oh_L^\times \times \hZ$, that is,
the profinite completion
of $L^\times$ \cite{32}. This extension is classified by an element of
\[
H^2(\Gal(L/\Q_p),L^\times) \cong \Z/n\Z \;,
\]
which also classifies those algebras, simple with center $\Q_p$, which
split after tensoring with
$L$. The extension in question generates this group, by a fundamental
result of local classfield
theory; but the division algebra with invariant $1/n$ also generates
this group, and the
associated group extension is the normalizer of $L^\times$ in $D^\times$.]
\end{subsubsect}

\begin{subsubsect}
\label{sec1.2.3}
The remarks in this subsection are a digression, but they
will be useful in \fullref{sec3.3}:
Since the normalizer of $L^\times$ acts as generalized automorphisms of
a Lubin--Tate group for
$L$, then for every $g \in \w(L^\ab/\Q_p)$ there is a power series
\[
[g](T) \in \oh_{L^\bnr}\llbracket T \rrbracket
\]
satisfying $[g_0]([g_1](T)) = [g_0 g_1](T)$, compatible with the natural
action of theWeil group
$\w(L^\ab/\Q_p)$ on the ring of integers $\oh_{L^\bnr}$ in $L^\bnr$. It
follows that $[g]( t) =
\alpha(g) T + \cdots$ defines a \emph{crossed} homomorphism
\[
\alpha \co \w(L^\ab/\Q_p) \to (L^\bnr)^\times \;,
\]
ie a map satisfying $\alpha(g_0 g_1) = \alpha(g_0) \cdot
\alpha(g_1)^{g_0}$, the superscript
denoting the action of the Weil group on $L^\bnr$ through its quotient
$\w(L^\nr/\Q_p)$.
When $L = \Q_q$ this implies the existence of an extension of the identity
homomorphism from
$\Q_q^\times$ to itself, to a crossed homomorphism from $\Q_q^\times
\ltimes \Z$ to
$(L^\bnr)^\times$. A corollary is the existence of a representation
$\Q_p^\bnr(1)$ of
$\w(\Q_q^\ab/\Q_p)$ on the completion $\Q_p^\bnr$ extending the action
of $W(\F_q)^\times$
by multiplication. [The completions in this construction are cumbersome,
and might be
unnecessary. Experts may know how to do without them, but I
don't.]
\end{subsubsect}
\end{subsect}

\begin{subsect}
\label{sec1.3}
The representations promised in the introduction are the
compositions
\[
\rho_q \co \w(\bQ_p/\Q_p) \to \w(\Q_q^\ab/\Q_p) \to D^\times
\]
with the second arrow coming from the Weil--Shafarevich theorem.

In the remainder of this section I will sketch a conjecture about the
relations between these
representations. The argument goes back to Serre's {\it Cohomologie
Galoisienne}, and I
believe that many people have thought along the lines below, but I don't
know of any place in
the literature where this is spelled out. It is based on a $p$--adic
analog of a conjecture of
Deligne, related to an older conjecture of Shafarevich (see Morava~\cite{34}; see also
Furusho~\cite{14}).

\begin{subsubsect}
We need some basic facts:
\begin{enumerate}
\item[(i)] in order that a pro--$p$--group be free, it is necessary and sufficient
that its
cohomological dimension be $\leq 1$ (see Serre~\cite[I.4.2,
Corollary~2]{38}), and

\item[(ii)] the maximal unramified extension of a local field with perfect
residue field has
cohomological dimension $\leq 1$ \cite[II.3.3, Ex~c]{38}.
\end{enumerate}

I will follow Serre's notation, which very similar to that used above. $K$
is a field complete
with  respect to a discrete valuation, with residue field $k$, eg a
finite extension of $\Q_p. \;
K^\nr$ will denote its maximal unramified extension, $K^s$ its (separable)
algebraic closure, and
$K^\tr$ will be the union of the tamely ramified Galois extensions
(ie with Galois group of order
prime to $p$) of $K$ in $K^s$; thus
\[
K^s \supset K^\tr \supset K^\nr \supset K \;.
\]
Eventually $K$ will be the quotient field $\Q_q$ of the ring $W(k)$
of Witt vectors for some
$k = \F_q$ with $q = p^f$ elements, but for the moment we can be more
general. We have
\[
\Gal(K^\tr/K^\nr) \cong \varprojlim \; \F_{p^n}^\times \cong \prod_{l
\neq p} \Z_l =
\hZ(1)(\neg p)
\]
\cite[II.4, Ex~2a]{38}; the term on the right denotes the component `away from
$p$' of the Tate
representation of $\Gal(\bk/k) \cong \hZ$, in which $1 \in \Z$ acts as
multiplication by $q = \#(k)$ \cite[II.5.6 Ex~1]{38}. The kernel in the
extension
\[
1 \to \Gal(K^s/K^\tr) = P \to \Gal(K^s/K^\nr) \to \Gal(K^\tr/K^\nr) \to 1
\]
is a pro--$p$--group \cite[Ex~2b]{38} closed in a group of cohomological
dimension one (by
assertion (ii)), hence itself of cohomological dimension one \cite[I.3.3,
Proposition~14]{38},
hence \emph{free} by assertion (i).

The extension in question is the \emph{inertia} group of $K$; it has a
natural $\Gal(\bk/k)$--action.
\end{subsubsect}

\begin{subsubsect}
From now on I will assume that $K = \Q_q$ is unramified over
$\Q_p$.

\begin{conj}
The extension displayed immediately above splits:
$\Gal(K^s/K^\nr) \cong P
\rtimes \hZ(1)(\neg p)$ is a semidirect product.
\end{conj}

[Since the kernel and quotient have relatively prime order, this would
be obvious if either were
finite. This may be known to the experts.]

Let $W(k)_0^\times = (1 + p W(k))^\times$ be the group of those units
of $W(k)$ congruent
to 1 mod $p$: the logarithm defines an exact sequence
\[
1 \to k^\times \to W(k)^\times \to W(k) \to k \to 0
\]
taking $W(k)_0^\times$ isomorphically to $pW(k)$. Let
\[
\W(\bk)^\times_0 = \varprojlim \; \{ W(k')_0^\times \;|\; k' \; {\rm
finite} \subset \bk \}
\]
be the limit under norms; it is a compactification of
$W(\bk)_0^\times$.

\begin{conj}
The topological abelianization of $P$ is naturally
isomorphic to $\W(\bk)_0^\times$.
\end{conj}

Alternately: Lazard's group ring
\[
\Z_p\llbracket P\rrbracket \cong \varprojlim \; \{ \hat{T}(W(k')_0^\times) \;|\: k' \;
{\rm finite} \; \subset \
\bk \}
\]
of $P$ of is isomorphic to the limit (under norms) of the system of
completed tensor algebras of
$pW(k') \cong W(k')_0^\times$. It is thus (hypothetically) a kind of
noncommutative Iwasawa
algebra \cite[I.1.5, Proposition~7, page~8]{38}.

If both these conjectures are true, then I can think of no natural way
for $\hZ(1)(\neg p)$
to act on $P$, so I will go the rest of the way and conjecture as well
that this action
is \emph{trivial}. I will also abbreviate the sum of these conjectures
as the assertion that
the absolute inertia group $\bI$ of $\Q_p$ is the product of $\hZ(1)(\neg
p)$ with the pro-free
pro--$p$--group generated by $W(\bk)_0^\times \cong pW(\bk)$.
\end{subsubsect}

\begin{subsubsect}
These conjectures seems to be compatible with other known
facts of classfield
theory. In particular, they would imply that
\[
\Gal(\bQ_p/\Q_q) \cong \bI \ltimes \hZ
\]
with $1 \in \Z$ acting as $\sigma(x) = x^q$ on $k$. This abelianizes to
$W(k)^\times \times
\hZ$, agreeing with Artin's reciprocity law, and if $q_1 = q_0^m$ this
would yield an exact
sequence
\[
1 \to \Gal(\bQ_p/\Q_{q_1}) \cong \bI \ltimes \hZ \to \Gal(\bQ_p/\Q_{q_0})
\cong \bI \ltimes
\hZ \to \Gal(\F_{q_1}/\F_{q_0}) \cong \Z/m\Z \to 0 \;.
\]
Finally, the composition
\[
\bI \to \bI_\ab \cong \W(\bk)_0^\times \times \hZ(1)(\neg p) \to
W(\F_q)^\times
\]
defines a compatible system of candidates for the quotient maps
\[
\rho_q \co \Gal(\bQ_p/\Q_p) \cong \bI \ltimes \hZ \to W(\F_q)^\times
\ltimes \hZ \cong
\Gal(\Q_q^\ab/\Q_p) \;.
\]
\end{subsubsect}

\section{Some more global representations}
\label{sec2}

This section is concerned with two very global objects, each of which
has some unfamiliar
features. The first subsection is concerned with the suspension spectrum
\[
S[BU] = \Sigma^\infty(BU_+)
\]
of the classifying space for the infinite unitary group, and its
interpretation (following Rognes~\cite{37}) as a
Hopf algebra object in the category of spectra. The second subsection
reviews some properties
of the Hopf algebra of quasisymmetric functions, following
Cartier~\cite{5} and Hoffman~\cite{23}. A
remarkable number of the
properties of the Hopf algebra of symmetric functions generalize to this
context (see Hazelwinkel~\cite{20}), and I will
try to keep this fact in focus. This Hopf algebra is conjectured to be
closely related to a certain
motivic group of interest in arithmetic geometry, and I have tried to
say a little about that, in
particular because it seems to overlap with recent work of Connes and
Marcolli, discussed in \fullref{sec2.3}.

\addtocounter{subsect}{1}
\setcounter{subsubsect}{0}
\begin{subsubsect}
Since $BU$ is an infinite loopspace, $S[BU]$ becomes an
$E_\infty$ ringspectrum
(or, in an alternate language, a commutative $S$--algebra); but $BU$
is also a {\it space}, with a
diagonal map, and this structure can be used to make this ringspectrum
into a Hopf algebra in
the category of $S$--modules.

In complex cobordism, the complete Chern class (see Adams~\cite{1})
\[
c_t = \sum c_I \otimes t^I = \phi^{-1} s_t \phi(1) \in MU^*(BU_+)
\otimes S_*
\]
(where $S_* = \Z[t_i \;|\; i \geq 1]$ is the Landweber--Novikov algebra,
with coaction $s_t$)
represents a morphism $S[BU] \to MU \wedge MU$ of ringspectra. Indeed,
algebra morphisms
from $S[BU]$ to the Thom spectrum $MU$ correspond to maps
\[
\C P^\infty_+ = MU(1) \to MU
\]
of spectra, and hence to elements of $MU^*(\C P^\infty_+)$. The coaction
\[
MU \cong S \wedge MU \to MU \wedge MU
\]
is also ring map, so the resulting composition
\[
S[BU]  \to MU = S \wedge MU \to MU \wedge MU
\]
is a morphism of ring spectra. On the other hand the Thom isomorphism
\[
\phi\co MU^*(BU_+) \cong MU^*(MU)
\]
satisfies
\[
\phi^{-1} s_t \phi(1) = c_t = \prod_i \e_i^{-1}t(\e_i)
\]
(where $t(\e) = \sum t_k \e^{k+1} \in MU^*(\C P^\infty_+) \otimes S_*$,
with $\e$ being the
Euler, or first Chern, class).

In fact this map is also a morphism of Hopf algebra objects:
\[
s_{t'}(\Delta c_t) = c_{t' \circ t} \otimes c_{t' \circ t} \in MU^*(BU_+
\wedge BU_+) \otimes
(S_* \otimes S_*) \; ,
\]
so the diagram
\[\xymatrix
{
{S[BU]} \ar[1,0]\ar[0,1]& {S[BU] \wedge S[BU]} \ar[0,1]& {(MU \wedge MU)
\wedge (MU
\wedge MU)} \ar[1,0]\\
{MU \wedge MU} \ar[0,2]& & {(MU \wedge MU) \wedge_{MU} (MU \wedge MU)}
}
\]
commutes.

We can thus think of $S[BU]$ as a kind of Galois group for the category
of $MU$--algebras
over $S$, or alternately for the category of complex-oriented
multiplicative cohomology
theories. In particular, if $E$ is an $MU$--algebra, algebra maps from
$S[BU]$ to $E$ define
elements of $E^0 (\C P^\infty_+)$; thus the group $\aut(E)$ of
multiplicative automorphisms of
$E$ maps to algebra homomorphisms from $S[BU]$ to $E$. The adjoint
construction thus
sends $S[BU]$ to the spectrum of maps from $\aut(E)$ (regarded naively,
as a set) to $E$; but
these maps can be regarded as $E_*$--valued functions on $\aut(E)$,
and hence as elements of
the coalgebra $E_*E$.

Note that $S[BU]$ is large, so the existence of an honest dual object
in the category of
spectra may be problematic. The remark above implies that this perhaps
nonexistent
group object admits the etale groupschemes $\oh_D^\times$ of
\fullref{sec1} as subgroups.
\end{subsubsect}

\begin{subsubsect}
The Landweber--Novikov algebra represents the group of
invertible power series
under composition, and it may be useful below to know that (since all
one-dimensional
formal groups over the rationals are equivalent) $(H_*(MU,\Q), H_*(MU
\wedge MU,\Q))$
represents the transformation groupoid defined by this group acting on
itself by
translation. On the other hand $H_*(BU,\Q)$ represents the group of formal
power series
with leading coefficient 1, under multiplication, and the induced morphism
\[
\spec (H_*(MU,\Q), H_*(MU \wedge MU,\Q)) \to \spec (H_*(S,\Q),H_*(BU,\Q))
\]
of groupoidschemes sends the pair $(g,h)$ of invertible series, viewed as
a morphism
from $h$ to $g \circ h$, to the translated derivative $g'(h(t))$. This is
essentially
the construction which assigns to a formal group, its canonical invariant
differential.
The chain rule ensures that this is a homomorphism: we have $(g,h) \circ 
(k, g \circ h) = (k \circ g, h)$, while
\[
g'(h(t)) \cdot k'((g \circ h)(t)) = (k \circ g)'(h(t)) \;.
\]
I'm indebted to Neil Strickland for pointing out the advantage of working
in this context
with $MUP$, that is, the spectrum $MU$ made periodic: its homotopy groups
represent the functor which
classifies formal group laws together with a coordinate, while $MUP \wedge
MUP$ represents
formal group laws with a \emph{pair} of coordinates. The unit classified by
$S[BU]$ is then
just the ratio of these coordinates.
\end{subsubsect}

\addtocounter{subsect}{1}
\setcounter{subsubsect}{0}
\begin{subsubsect}
The $\Z$--algebra of symmetric functions manifests itself in
topology as the integral
homology of $BU$; it is a commutative and cocommutative Hopf algebra
(with a canonical
nondegenerate inner product \cite{20} which, from the topological point of
view, looks quite
mysterious). Over $\Q$ it is the universal enveloping algebra of an
abelian graded Lie algebra
with one generator in each even degree.

The graded ring $\QSymm_*$ of {\it quasi} symmetric functions is the
commutative Hopf
algebra dual to the free associative algebra on
(noncommutative!) generators $Z_k$ of degree
$2k$, with coproduct
\[
\Delta Z_k = \sum_{i+j=k} Z_i \otimes Z_k \;;
\]
it is thus the universal enveloping algebra for a free graded Lie algebra
$\ff_\Z$, with one
generator in each (even) degree (see Hazelwinkel~\cite{19}). There is a natural monomorphism
embedding the
symmetric functions in the quasisymmetric functions (see
Cartier~\cite[Section~2.4]{5}), dual to
the map on
enveloping algebras defined by the homomorphism from the free Lie algebra
to its
abelianization. This then defines a morphism
\[
S[BU] \to H\Z \otimes\QSymm_*
\]
of (Hopf) ringspectra.

[In \cite{3a}, which appeared after this paper was
written, Baker and Richter show that the ringspectrum $S[\Omega \Sigma \C
P^\infty_+]$ has many properties one might expect of a dual (but see the
cautionary remarks in the paragraphs above) to the algebra of functions on
$\Gal_\mot$.]
\end{subsubsect}

The group-valued functors represented by such Hopf algebras are very
interesting, and have a
large literature. In what follows I will simplify by tensoring everything
with $\Q$, which will be
general enough for anything I have to say. The formal (Magnus) completion
of the $\Q$--algebra
of noncommuting power series is dual to the algebra of functions on a
free prounipotent (see Deligne~\cite[Section~9]{10}) groupscheme $\fF$. In this context, a grading on a Lie algebra can
be reinterpreted as an
action of the multiplicative groupscheme $\G_m$, which sends an element
$x$ of degree $d$ to
$\lambda^d x$, where $\lambda$ is a unit in whatever ring we're working
with, so we can
describe the map above as defining a morphism
\[
\fF \ltimes \G_m \to \Gal_\hot
\]
of group objects of some sort, over $\Q$.

\begin{subsubsect}
In arithmetic geometry there is currently great interest in
a groupscheme
\[
\Gal_\mot = \fF_\odd \ltimes \G_m
\]
which is conjectured to be (isomorphic to) the motivic Galois group of
a certain Tannakian
category, that of mixed Tate motives over $\Z$ (see Deligne
and Goncharov~\cite{10,12}). Usually $\fF_\odd$ is taken to be the
prounipotent group defined by the free graded $\Q$--Lie algebra $\ff_\odd$
on generators of
degree $4k+2$ (hence `odd' according to the algebraists' conventions),
with $k \geq 1$; but
there are reasons to allow $k=0$ as well. We can regard $\fF_\odd$
as a subgroupscheme of
$\fF$, by regarding $\ff_\odd$ as a subalgebra of $\ff$.

A Tannakian category is, roughly, a suitably small $k$--linear abelian
category with tensor
product and duality -- such as the category of finite-dimensional linear
representations of a
proalgebraic group over a field $k$. Indeed, the main theorem of the
subject (see Deligne~\cite{11}) asserts that
(when $k$ has characteristic zero) any Tannakian category is of this form;
then the relevant
group is called the motivic group of the category. Present technology
extracts mixed Tate
motives from a certain triangulated category of (pieces of) algebraic
varieties, constructed as a
subcategory of the tensor triangulated category of more general motives
(see Deligne and Goncharov~\cite{12} and Voevodsky~\cite{44}).

There are more details in \fullref{sec3} below, but one of the points of this
paper is that the language of
such tensor categories can be quite useful in more general
circumstances. For example, the
category of complex-oriented multiplicative cohomology theories behaves
very much like (a
derived category of) representations, with $S[BU]$ as its motivic
group. Away from the prime
two, the fibration
\[
SO/SU_+  \to BSU_+ \to BSO_+
\]
(defined by the forgetful map from $\C$ to $\R$) splits (even as maps
of infinite loopspaces).
Over the rationals, this is almost trivial: it corresponds to the
splitting of the graded abelian Lie
algebra with one generator in each even degree, into a sum of two such
Lie algebras, with
generators concentrated in degrees congruent to 0 and 2 mod 4
respectively. The composition
(the first arrow is the projection of the $H$--space splitting, and the
second corresponds to the
abelianization of a graded free Lie algebra).
\[
S[BU] \to S[SO/SU] \wedge S[\C P^\infty] \to H\Q \otimes \QSymm_\odd
\]
is the candidate, promised in the introduction, for a natural
representation
\[
\Gal_\mot \to \Gal_\hot
\]
(over $\Q$, of course!).
\end{subsubsect}

\begin{subsubsect}
The degree two $( = 2 (2 \cdot 0 + 1))$ generator in
$\fF_\odd$ is closely related to
the $S[\C P^\infty]$ factor in the ring decomposition above; both
correspond to exceptional
cases that deserve some explanation.
\end{subsubsect}

The symmetric \emph{functions} are defined formally as an inverse limit
of rings of symmetric \emph{polynomials} in finitely many variables $x_n, \; n \geq 1$. There are
thus interesting maps from
the symmetric functions to other rings, defined by assigning interesting
values to the $x_n$; but
because infinitely many variables are involved, issues of convergence
can arise. For example: if
we send $x_n$ to $1/n$, the $n$th power sum $p_n$ maps to $\zeta(n)
\in \R$ \ldots\ as long as
$n >1$, for $s = 1$ is a pole of $\zeta(s)$.

This is a delicate matter (see Cartier~\cite[Section~2.7]{5} and
Hoffman~\cite{23}), and it turns out to be very
natural to send $p_1$ to
Euler's constant $\gamma$. In fact this homomorphism extends, to define
a homomorphism
from $\QSymm$ to $\R$, whose image is sometimes called the ring $\MZN_*$
of `multizeta
numbers'; it has a natural grading. It is classical that for any positive
integer $n$,
\[
\zeta(2n) = - \half B_{2n} \; \gamma_{2n}(2 \pi i) \in \Q(\pi) \;,
\]
(where $\gamma_k$ denotes the $k$th divided power) and in some contexts
it is natural to work
with the even-odd graded subring of $\C$ obtained by adjoining an
invertible element $(2\pi
i)^\pm$ to $\MZN_*$. It is known (see Hain~\cite[Section~4]{18}) that the multizeta numbers
are periods of algebraic
integrals, and that the ring generated by all such periods is a Hopf
algebra, closely related to the
algebra of functions on the (strictly speaking, still hypothetical)
motivic group of all motives over
$\Q$ (see Kontsevich~\cite{28}); but it is thought that Euler's constant is probably not a
period. Nevertheless, from the
homotopy-theoretic point of view presented here, it appears quite
naturally.
\end{subsect}

\addtocounter{subsect}{1}
\setcounter{subsubsect}{0}

\begin{subsubsect}
These multizeta numbers may play some universal role in the general theory
of asymptotic expansions; in any case, they appear systematically in
renormalization theory. Connes and Marcolli, building on earlier work
of Connes and Kreimer~\cite{7a,7b}, Broadhurst, and others, have put
this in a Galois-theoretic framework. This subsection summarizes some
of their work.

Classical techniques [of Bogoliubov, Parasiuk, Hepp, and Zimmerman]
in physics have achieved an impressive level of internal consistency,
but to mathematicians they lack conceptual coherence.  Starting with
a suitable Lagrangian density, Connes and Kreimer define a graded
Lie algebra $\fg_*$ generated by a class of Feynman graphs naturally
associated to the interactions encoded by the Lagrangian. They interpret
the BPHZ dimensional regularization procedure as a Birkhoff decomposition
for loops in the associated prounipotent Lie group $\cG$, and construct
a universal representation of this group in the formal automorphisms of
the line at the origin, yielding a formula for a reparametrized coupling
constant which eliminates the divergences in the (perturbative) theory
defined by the original Lagrangian.

Dimensional regularization involves regarding the number $d$ of space-time
dimensions as a special value of a complex parameter; divergences are
interpreted as poles at its physically significant value. Renormalization
is thus expressed in terms much like the extraction of a residue,
involving a simple closed curve encircling the relevant value of $d$ (the
source of the loop in $\cG$ mentioned above). Connes and Marcolli~\cite{8}
reformulate such data in geometric terms, involving flat connections on a
$\G_m$--equivariant $\cG$--bundle over a certain `metaphysical' (not their
terminology!) base space $B$. From this they define a Tannakian category
(of flat, equisingular vector bundles with connection over $B$), and
identify its motivic group as $\fF \ltimes \G_m$. Their constructions
define a representation of the motivic group in $\cG$, and hence in the
group of reparametrizations of the coupling constant.
\end{subsubsect}

\begin{subsubsect}
Besides the complex deformation of the space-time dimension,
the base space $B$
encodes information about the mass scale. It is a (trivial, but not
naturally trivialized, see the end of \cite[Section~2.13]{9})
principal bundle
\[
\G_m \to B \to \Delta
\]
over a complex disk $\Delta$ centered around the physical dimension $d
\in \C$. [This disk is
treated as infinitesimal but there may be some use in thinking of it as
the complement of infinity in
$\C P_1$.] The renormalization group equations \cite[Section~2.9]{8} are reformulated
in terms of the
(mass-rescaling) $\G_m$--action on the principal bundle $\cG \times B$
(the grading on $\fg_*$
endows $\cG$ with a natural $\G_m$--action), leading to the existence of
a unique gauge-equivalence
class of flat $\G_m$--equivariant connection forms $\lambda \in
\Omega^1(B,\fg)$,
which are equisingular in the sense that their restrictions to sections
of $B$ (regarded as a
principal bundle over $\Delta$) which agree at $0 \in \Delta$ are mutually
(gauge) equivalent.

A key result \cite[Theorem~2.25, Section~2.13]{8} characterizes such forms $\lambda$
in terms of a graded
element $\beta_* \in \fg_*$ corresponding to the beta-function of
renormalization group theory.
In local coordinates ($z \in \Delta$ near the basepoint, $u \in \G_m$)
we can write
\[
\lambda(z,u) = \lambda_0(z,u) \cdot dz + \lambda_1(z,u) \cdot u^{-1}du
\]
with coefficients $\lambda_i \in \C\{u,z\}[z^{-1}]$ allowed singularities
at $z = 0$; but flatness
and equivariance imply that these coefficients determine each other. The
former condition \cite[Equation~2.166]{8}
can be stated as
\[
\partial_z \lambda_1 = H\lambda_0 - [\lambda_0,\lambda_1]
\]
where the grading operator
\[
H = (u \partial_u)|_{u=1}
\]
is the infinitesimal generator of the $\G_m$--action; it sends $u^k$
to $ku^k$. Regularity of
$\lambda_0$ at $u=0$ implies that
\[
\lambda_0 = H^{-1}[\partial_z \lambda_1 + [H^{-1} \partial_z
\lambda_1,\lambda_1] + \cdots]
\]
is determined, at least formally, by knowledge of $\lambda_1$ in the
fiber direction. Connes and
Marcolli's solution \cite[Theorems~2.15 and~2.18, Equation~2.173]{8} of the
renormalization group equations
imply that
\[
\lambda_1(u,z) = - z^{-1}u^*(\beta)
\]
for some unique $\beta \in \fg$. This characterizes a universal (formal)
flat equisingular
connection $\lambda(\beta)$ on $B$, related to the universal singular
frame of \cite[Section~2.14]{8}.
\end{subsubsect}

\begin{subsubsect}
The Lie algebra $\bV$ of the group of formal diffeomorphisms
of the line at the origin
has canonical generators $v_k = u^{k+1}\partial_u$ satisfying $[v_k,v_l] =
(l - k) v_{k+1}, \; k,l
\geq 1$, so a graded module with an action of such operators defines a
flat equisingular vector
bundle over $B$ with $\beta_* = v_*$. A commutative ringspectrum $E$
with $S[BU]$--action
defines a multiplicative complex-oriented cohomology theory, and in
particular possesses an
$MU$--module structure. On rationalized homotopy groups the morphism
\[
E \to E \wedge S \to E \wedge MU
\]
defines an $S_*$--comodule structure map
\[
E_\Q^* \to E^*_\Q \otimes_{MU^*_\Q} MU^*MU_\Q = E^*_\Q \otimes S_*
\]
and thus an action of the group of formal diffeomorphisms; differentiating
this action assigns $E^*_\Q$ an action of $\bV$. The smash product of
two complex-oriented spectra is another such thing, and the functor from
$S[BU]$--representations to $\bV$ representations takes this product
to the usual tensor product of Lie algebra representations. This,
together with the universal connection constructed above, defines a
monoidal functor from $S[BU]$--representations to the Tannakian category
of flat equisingular vector bundles over $B$ of \cite[Section~2.16]{8},
and thus a (rational) representation of its motivic group in $S[BU]$.
\end{subsubsect}

\begin{subsubsect}
\label{sec2.3}
Connes and Marcolli note that their motivic group is isomorphic to the
motivic group for mixed Tate motives over the Gaussian integers (which
has generators in all degrees, unlike that for mixed Tate motives over
the rationals), but they do not try to make this isomorphism canonical.
It is striking to me that their formulas \cite[Equation~2.137]{8} actually
take values in the group of \emph{odd} formal diffeomorphisms. Presumably
$\Gal(\Q(i)/\Q)$ acts naturally on the motivic group of mixed Tate
motives over the Gaussian integers, and it seems conceivable that these
constructions might actually yield a representation of $\fF_\odd \ltimes
\G_m$ in the group of odd formal diffeomorphisms \cite{34}; but I have
no hard evidence for this.

Note that the function field $\Q(\Gamma(z),z)$ admits the endomorphism
\[ \tau(z) = z + 1 \]
making it into a difference field, with a flow defined by
\[ \exp(t \tau) \; \Gamma(z) = (1 - t)^{-z} \Gamma(s) \]
satisfying
\[ \exp(t \tau) \; z = (1 - t)^{-z} \;; \]
as well -- but beware, this group does not act multiplicatively! Taylor
expansion at the origin defines a morphism
\[ \Gamma(z)^{-1} \mapsto z \; \exp\biggl(- \sum_{k \geq 1}
\frac{\zeta(k)}{k} z^k\biggr) \co \Q(\Gamma(z),z) \to \R((z)) \]
of fields, compatible with this flow; this can be interpreted as a
Tannakian fiber functor, from some category of vector bundles over
the sphere endowed with a shift operator. Its motivic group is closely
related to the automorphisms of the transcendental extension of $\Q$
generated by the odd zeta-values.

Dimensional regularization replaces certain divergent Feynman integrals
with expressions involving Gamma-functions. This `asymptotic point'
sends these function-field expressions to formal power series with
coefficients in the field generated over the rationals by zeta-values.

\end{subsubsect}

\section{Towards fundamental groupoids of tensor triangulated categories}
\label{sec3}

\addtocounter{subsect}{-1}
\begin{subsect}
A Tannakian category $\cA$ is a $k$--linear abelian category,
where $k$ is a field,
possessing a coherently associative and commutative tensor product
$\otimes$ (see Deligne~\cite[Section~2.5]{11});
moreover, $\cA$ should be small enough: its objects should be of finite
length, its Hom-objects
should be finite-dimensional over $k$, the endomorphism ring of the
identity object for the tensor
product should be $k$, and it should admit a good internal duality
\cite[Section~2.12]{12}. The specifically
Tannakian data, however, consists of a nontrivial exact $k$--linear functor
\[
\omega \co \cA \to (k-\Vect)
\]
which is monoidal in the sense that
\[
\omega(X \otimes Y) \cong \omega(X) \otimes_k \omega(Y) \;.
\]
Let $\aut_\omega (k)$ be the group of multiplicative automorphisms of
$\omega$, ie of natural
transformations of $k$--module-valued functors from $\omega$ to itself,
which commute with
these multiplicativity isomorphisms. More generally, the multiplicative
automorphisms
\[
A \mapsto \aut_\omega (A) = \aut(\omega \otimes_k A)
\]
of $\omega \otimes_k A$ define a group-valued functor on the category
of commutative
$k$--algebras.

It sometimes happens (eg when $k$ is of characteristic zero) that this
functor is \emph{representable} by a suitable Hopf algebra $\cH$, ie
\[
\aut(\omega \otimes_k A) \cong \Hom_{k - {\rm alg}}(\cH,A) \;;
\]
then $\omega$ lifts to a functor from $\cA$ to the category of
finite-dimensional representations
of the affine (pro)algebraic group represented by $\cH$, and it may be
that we can use this lift to
identify $\cA$ with such a category of representations
\cite[Section~7]{11}. In any
case, when the
$k$--groupscheme
\[
\spec \cH = \pi_1(\spec \cA,\omega)
\]
exists, it is natural to think of it as a (`motivic') kind of fundamental
group for the Tannakian
category, with the `fiber functor' $\omega$ playing the role of
basepoint.
\end{subsect}

\begin{subsect}
Here are some illustrative examples, and variations on this
theme:

\begin{subsubsect}
The category of local systems (ie of finite-dimensional flat
$k$--vector spaces) over a
connected, suitably locally connected topological space $X$ is Tannakian:
a basepoint $x \in X$
defines an exact functor, which sends the system to its fiber over $x$,
and we recover the
fundamental group $\pi_1(X,x)$ (or, more precisely, its `envelope'
or best approximation by a
proalgebraic group over $k$) as its automorphism object (see
Toen~\cite{42}).
\end{subsubsect}

\begin{subsubsect}
The functor which assigns to a finite-dimensional $\Q$--vector
space $V$, the
$p$--adic
vector space $V \otimes_\Q \Q_p$, defines a Tannakian structure on
$\Q-\Vect$, with
$\Gal(\bQ_p/\Q_p)$, regarded as a profinite groupscheme over $\Q_p$,
as its motivic
fundamental group.
\end{subsubsect}

\begin{subsubsect}
One possible variation involves fiber functors taking values
in categories more general
than vector spaces over a field. The monoidal category of even-odd graded
(or `super') vector
spaces is an important example. When this works, one gets a Hopf algebra
object (corresponding
to a `super' groupscheme) in the enriched category; the motivic group
in this case is the
multiplicative group $\mu_2$ of square roots of unity (suitably
interpreted \cite[Section~8.19]{11}).
\end{subsubsect}

\begin{subsubsect}
\label{sec3.1.4}
The automorphisms of cohomology with $\F_p$ coefficients,
which is defined not on
an abelian category but on a tensor triangulated one, eg finite spectra,
is a more exotic and
perhaps more compelling example: the even-degree cohomology is a
representation of the
groupscheme defined by the functor
\[
A \mapsto \{ \sum_{k \geq 0} a_k T^{p^k} = a(T) \in A\llbracket
T\rrbracket \; | \:
a_0 \in A^\times \} \;;
\]
its Hopf algebra is dual to the algebra defined by Steenrod's reduced
$p$th powers.

The group operation here is composition of power series. The Hopf algebra
dual to the full
Steenrod algebra can be recovered by working with $\mu_2$--graded (super)
vector spaces.
Note that the groupscheme defined above contains the `torus'
\[
A \mapsto A^\times = \Hom_{k-{\rm alg}}(k[t_0^{\pm1}],A)
\]
as the subgroup of series of the form $a_0T, \; a_0 \in A^\times$;
this action by the
multiplicative groupscheme $\G_m$ allows us to recover the (even part of)
the grading on the
cohomology in `intrinsic' terms. However, we're not as lucky here as in
the preceding case: the
stable homotopy category (at $p$) is not the same as the category of
modules over the Steenrod
algebra. The Adams spectral sequence tells us that we have to take higher
extensions into
consideration.
\end{subsubsect}

\begin{subsubsect}
These motivic fundamental groups are in a natural sense
functorial \cite[Section~8.15]{11}:
Given a commutative diagram
\[
\xymatrix{
{\cA_0}\ar[r]^{\eta}\ar[d]^{\omega_0}&{\cA_1}\ar[d]^{\omega_1}\\
{(k_0-\Vect)}\ar[r]^{\tilde{\eta}}&{(k_1-\Vect)}}
\]
of Tannakian categories $\cA_i$, fiber functors $\omega_i$, and fields
$k_0 \to k_1$, with $\eta$
exact $k_0$--linear multiplicative and $\tilde{\eta} = - \otimes_{k_0}
k_1$, there is a natural
homomorphism
\[
\eta^* \co \pi_1(\spec \cA_1,\omega_1) \to \pi_1(\spec \cA_0,\omega_0)
\times_{k_0} k_1
\]
of groupschemes over $k_1$ constructed as follows: if $\alpha \co \omega_1
\cong \omega_1$ is a
multiplicative automorphism of $\omega_1$ (perhaps after some base
extension which I won't
record), then
\[
\alpha \circ \eta \co \tilde{\eta} \circ \omega_0  = \omega_1 \circ \eta
\cong  \omega_1 \circ \eta =
\tilde{\eta} \circ \omega_0
\]
is an element of $\pi_1(\spec \cA_0,\omega_0)$, pulled back by
$\tilde{\eta}$ to become a
groupscheme over $k_1$.

As \fullref{sec3.1.4} suggests, it is tempting to push these constructions in
various directions; in
particular, if $\cT$ is a suitably small tensor triangulated category
(ie a triangulated category with
tensor product satisfying reasonable axioms (for example that tensoring
with a suitable object
$\Sigma$ represents suspension) \cite{25}), and $\omega$ is a multiplicative
homological functor (taking
distinguished triangles to long exact sequences) with values in $k$--vector
spaces, we can consider
its multiplicative automorphisms, as in the Tannakian case; and we might
hope that if $\cT$ is
some kind of derived category of $\cA$, then we could try to reconstruct
$\cT$ in terms of a derived category of representations of the
automorphism group of the homological fiber functor $\omega$ (see Franke
\cite{14}, Neeman \cite{38}).
(If we ask that the automorphisms behave reasonably under
suspension, then the
automorphism group will contain a torus encoding gradings, as in
\fullref{sec3.1.4}.)
\end{subsubsect}
\end{subsect}

\begin{subsect}
What seems to be missing from this picture is a compatible understanding
of $\pi_0$. If $\fp$ is a prime ideal in a commutative noetherian ring $A$
then $A/\fp$ is a domain with quotient field $Q(A/\fp)$, and the composite
\[
\omega_\fp \co A \to A/\fp \to Q(A/\fp) = k(\fp)
\]
determines $\fp$, so the prime ideal spectrum $\spec A$ is a set of
equivalence classes of
homomorphisms from $A$ to (varying) fields $k$. Such a homomorphism
lifts to define a
triangulated functor $\otimes^L_A k(\fp)$ from the derived category of
$A$--modules to the
derived category of vector spaces over $k(\fp)$, taking (derived) tensor
products to (derived)
tensor products, and (as in the case of ordinary cohomology above)
we can consider the functor
of automorphisms associated to such a generalized point. $G$--equivariant
stable homotopy
theory is another example of a category with a plentiful supply of natural
`points', corresponding
to conjugacy classes of closed subgoups of $G$.

It would be very useful if we could construct, for suitable $\cT$,
a fundamental groupoid $\bpi (\spec \cT)$ with objects corresponding
to equivalence classes of multiplicative homological fiber functors
$\omega$, and morphisms coming from $\pi_1(\spec \cT,\omega)$. The
lattice of Bousfield localizations (see Hovey and Palmieri~\cite{24})
of a triangulated category is in some ways analogous to the Boolean
algebra of subsets of the spectrum of an abelian category (see
Balmer~\cite{3}, Hopkins~\cite{21}, Krause~\cite{29} and Neeman~\cite{36}) but we do not
seem to understand, in any generality, how to identify in it a sublattice
corresponding to the open sets of a reasonable topology. Moreover,
there seem to be subtle finiteness issues surrounding this question:
constructing a good $\pi_1$ at $\omega$ may require identifying
a suitable subcategory of $\omega$--finite or coherent objects (see
Christensen--Hovey~\cite{6} and Hovey--Strickland~\cite[Section~8.6]{26}), and we might hope that
a reasonable topology on $\spec \cT$ would suggest a natural theory of
adeles (see Kapranov~\cite{27}) associated to chains of specializations.
\end{subsect}

\begin{subsect}
\label{sec3.3}
This is all quite vague; perhaps some examples will be
helpful.

The structure map $S \to H\Z$ defines a monoidal pullback functor
\[
- \wedge H\Z \co {\rm (Spectra)} \to D(\Z - {\rm Mod}) \;;
\]
simply regarding a $\Z$--algebra as an $S$--algebra is not monoidal:
$H(A \otimes B)$ is not the
same as $HA \wedge HB$.

\begin{subsubsect}
Let $\hK^*_p( - ) = K^*( - ) \otimes \Z_p$ denote classical
complex $K$--theory,
regarded as a functor on finite spectra and $p$--adically completed. The
Adams operations,
suitably normalized, define an action of $\Z^\times_p$ by stable
multiplicative automorphisms (see Sullivan~\cite{40}),
and the Chern character isomorphism
\[
\hK^*_p( - ) \otimes \Q \to H^*( - , \Q_p)
\]
identifies its eigenspaces with the graded components defined by ordinary
cohomology.

This can be expressed in the language developed above, as follows:
the fiber functor $\omega_0 =
\hK_p \otimes \bQ_p$ on finite spectra is ordinary cohomology in disguise,
and so has
automorphism group $\G_m \ltimes \Gal(\bQ_p/\Q_p)$, as does
\[
\omega_1 \co M \mapsto H^+(M \otimes \bQ_p) \co D(\Z-{\rm Mod}) \to
(\bQ_p-\Vect) \;.
\]
The induced morphism
\[
\G_m \ltimes \Gal(\bQ_p/\Q_p) \to \G_m \ltimes \Gal(\bQ_p/\Q_p)
\]
of motivic fundamental groups sends $(u,g) \in \G_m \ltimes
\Gal(\bQ_p/\Q_p)$ to $(g_\ab u,
g)$,
where
\[
g \mapsto g_\ab\co \Gal(\bQ_p/\Q_p) \to \Gal(\bQ_p/\Q_p)_\ab \to \Z^\times_p
\]
is the $p$--adic cyclotomic character; in Hopf algebra terms, this
corresponds to the
homomorphism
\[
u \mapsto 1 - t \co \Q_p[u^{\pm 1}] \to \Z_p\llbracket t\rrbracket[\mu_{p-1}] \otimes \Q
\cong \Z_p\llbracket\Z_p^\times\rrbracket
\otimes \Q
\]
of Iwasawa theory.
\end{subsubsect}

\begin{subsubsect}
This has a conjectural generalization to the spectra obtained
by specializing $E_n$,
with
\[
E_n^*(S) = \Z_p\llbracket v_1,\ldots,v_{n-1}\rrbracket[u^{\pm 1}]
\]
at $v_i \mapsto 0, \; 1 \leq i \leq n-1$. [These spectra correspond to
the Honda formal group of \fullref{sec1.1.1}; they have some kind (see
Goerss and Hopkins~\cite{15}) of multiplicative structure. It is
convenient to call the associated
cohomology functors $K(n)^*( - ;\Z_p)$.]

It is natural \cite{33} to expect that the group of multiplicative
automorphisms of this functor maps,
under reduction modulo $p$, to the normalizer of the maximal torus
$\Q_q^\times$ in
$D^\times$, as in \fullref{sec1.2.1}; in other words, precisely those automorphisms
of $K(n)^*( - ;\F_p)$
which lift to automorphisms of the formal group law of $K(n)^*( -;\Z_p)$
lift to multiplicative
automorphisms of the whole functor. It would follow from this, that
\[
K(n)^*( -;\Z_p) \otimes_{\Z_p} \Q_p^\bnr \cong \oplus_{k \in \Z} H^*( - ;
\Q_p^\bnr(1)^{\otimes k}) \;,
\]
where $\Q_p^\bnr(1)$ is the representation of $\Gal(\Q_q^\ab/\Q_p)$
(and thus of
$\Gal(\bQ_p/\Q_p)$) constructed in \fullref{sec1.2.3}.  The results of the preceding
section would then
generalize, with the role of the cyclotomic character now being played
by the crossed
homomorphism $\alpha$.
\end{subsubsect}
\end{subsect}

\addtocounter{subsect}{1}
\setcounter{subsubsect}{0}
\begin{subsubsect}
This suggests that the remarks above about $\pi_0$ might be
slightly naive: $\pi_0$
might make better sense as a small diagram category than as a topological
space.
Equivariant stable homotopy theory provides further evidence for this,
as do derived categories
of representations of quivers in groups.

The conjecture above asserts that the functor
\[
A \mapsto \aut(K(n)^*( -; \Z_p) \otimes_{\Z_p} A) \co \text{ (\emph{flat}
$A$--{\rm Mod})} \to \text{({\rm Groups})}
\]
is (i) representable (ii) by $\Gal(\Q^\ab_q/\Q_p)$. Enlarging the useable
class of fiber functors
in this way, allowing values in (say) modules over discrete valuation
rings, would lead to the
existence of specialization homomorphisms such as
\[
\pi_1(\spec \Z_p) \to \pi_1(\spec S,K(n)^*( - ; \Z_p)) \to \pi_1(\spec
S,K(n)^*( - ;\F_p)) \;.
\]
\end{subsubsect}

\begin{subsubsect}
This paper is so much a wish list, that while we're at it we might as well
ask for a way to associate automorphism groups to suitable open subsets
or subdiagrams of our hypothetical $\pi_0$, together with homomorphisms
between them defined by inclusions.  One of the beauties of Grothendieck's
original account of the fundamental group is his theory of specialization
\cite[Section~X.2.3]{16} which, as he notes, has no immediate analog in
pure algebraic topology. There is evidence (see Ando, Morava and Sadofsky
\cite{2} and Torii \cite[Section~4]{43})
that some generalization of his theory to a derived category context
could accomodate a very general form of the theory of vanishing cycles:

Classically $i\co Y \to X$ is the inclusion of a closed subscheme, and
$j\co X - Y \to X$ is the inclusion of its complement; for example $X$
might be an object over the spectrum of a discrete valuation ring, with
$Y$ the fiber above the closed point, and $X - Y$ the fiber above the
generic point. In terms of derived categories and the various functors
relating them, Bousfield localization along $Y$ is the functor
\[ L_Y = j_*j^*\co D(X) \to D(X) \;, \]
see B\"okstedt--Neeman \cite[Section~6]{4} and Hopkins--Gross \cite{22};
it is related to the Grothendieck local cohomology functor $i_*i^!$
by an exact triangle
\[ \dots \to i_*j^! \to {\rm Id} \to j_*j^* \to \dots \]
[since we're working in derived categories, I won't bother to indicate
that all the functors have been suitably (left or right) derived].

The vanishing cycle functor (see Grothendieck \cite[Section~2.2]{17a,17b}
and Consani \cite{9}) is the related composition
\[ i^*j_*\co D(X - Y) \to D(Y) \;; \]
although for our purposes it might be better to think of this as the 
\emph{nearby} cycles
functor. The two constructions are related by the commutative diagram
\[ \xymatrix { D(X) \ar[r]^{j_*j^*} \ar[d]^{j^*} & D(X) \ar[d]^{i^*}\\
D(X - Y) \ar[r]^{i^*j_*} & D(Y) \;.  } \]
in other words the nearby cycles functor is roughly Bousfield localization
restricted to the open stratum. It would be nice to have an equivariant
version of this, or (in another language) a form adapted to stacks \cite{36}.
\end{subsubsect}

\bibliographystyle{gtart}
\bibliography{link}

\begin{thebibliography}{}
\providecommand\bibmarginpar{\leavevmode\marginpar}
\def\urlstyle#1{{\tt #1}}

\bibitem{1}
\textbf{J\,F Adams}, \emph{Stable homotopy and generalised homology}, Chicago
  Lectures in Mathematics, University of Chicago Press, Chicago, IL (1995)
  \xox{MR}{1324104}

\bibitem{2}
\textbf{M Ando}, \textbf{J Morava}, \textbf{H Sadofsky},
  \href{http://dx.doi.org/10.2140/gt.1998.2.145} {\emph{Completions of
  $\mathbb{Z}/(p)$--{T}ate cohomology of periodic spectra}}, Geom. Topol. 2
  (1998) 145--174 \xox{MR}{1638030}

\bibitem{3a}
\textbf{A Baker}, \textbf{B Richter}, \emph{Quasisymmetric functions from a
  topological point of view}  \xox{arXiv}{math.AT/0605743}

\bibitem{3}
\textbf{P Balmer}, \emph{The spectrum of prime ideals in tensor triangulated
  categories} \xox{arXiv}{math.CT/0409360}

\bibitem{4}
\textbf{M B{\"o}kstedt}, \textbf{A Neeman}, \emph{Homotopy limits in
  triangulated categories}, Compositio Math. 86 (1993) 209--234
  \xox{MR}{1214458}

\bibitem{5}
\textbf{P Cartier}, \emph{Fonctions polylogarithmes, nombres polyz\^etas et
  groupes pro-unipotents}, Ast\'erisque  (2002) Exp. No. 885, viii, 137--173
  \xox{MR}{1975178}

\bibitem{6}
\textbf{J\,D Christensen}, \textbf{M Hovey},
  \href{http://muse.jhu.edu/journals/american_journal_of_mathematics/v122/122.%
2christensen.pdf} {\emph{Phantom maps and chromatic phantom maps}}, Amer. J.
  Math. 122 (2000) 275--293 \xox{MR}{1749049}

\bibitem{7a}
\textbf{A Connes}, \textbf{D Kreimer}, \emph{Renormalization in quantum field
  theory and the Riemann--Hilbert problem I: the Hopf algebra structure of
  graphs and the main theorem} \xox{arXiv}{hep-th/9912092}

\bibitem{7b}
\textbf{A Connes}, \textbf{D Kreimer}, \emph{Renormalization in quantum field
  theory and the Riemann--Hilbert problem II: the $\beta$--function,
  diffeomorphisms and the renormalization group} \xox{arXiv}{hep-th/9912092}

\bibitem{8}
\textbf{A Connes}, \textbf{M Marcolli}, \emph{From physics to number theory via
  noncommutative geometry II: renormalization, the Riemann--Hilbert
  correspondence, and motivic Galois theory} \xox{arXiv}{hep-th/0411114}

\bibitem{9}
\textbf{C Consani}, \href{http://dx.doi.org/10.1023/A:1000362027455}
  {\emph{Double complexes and {E}uler $L$--factors}}, Compositio Math. 111
  (1998) 323--358 \xox{MR}{1617133}

\bibitem{10}
\textbf{P Deligne}, \emph{Le groupe fondamental de la droite projective moins
  trois points}, from: ``Galois groups over $\mathbb{Q}$ (Berkeley, CA,
  1987)'', Math. Sci. Res. Inst. Publ. 16, Springer, New York (1989)  79--297
  \xox{MR}{1012168}

\bibitem{11}
\textbf{P Deligne}, \emph{Cat\'egories tannakiennes}, from: ``The Grothendieck
  Festschrift, Vol.\ II'', Progr. Math. 87, Birkh\"auser, Boston (1990)
  111--195 \xox{MR}{1106898}

\bibitem{12}
\textbf{P Deligne}, \textbf{A Goncharov}, \emph{Groupes fondementaux motiviques
  de Tate mixte} \xox{arXiv}{math.NT/0302267}

\bibitem{14}
\textbf{H Furusho}, \href{http://dx.doi.org/10.1007/s00222-003-0320-9}
  {\emph{$p$--adic multiple zeta values. {I}. $p$--adic multiple polylogarithms
  and the $p$--adic {KZ} equation}}, Invent. Math. 155 (2004) 253--286
  \xox{MR}{2031428}

\bibitem{15}
\textbf{P Goerss}, \textbf{M Hopkins}, \emph{Moduli spaces for structured ring
  spectra}, Hopf preprint.
\ Available at \setbox0\hbox{\makeatletter\@url
{http://hopf.math.purdue.edu/cgi-bin/generate?/\penalty -100\unskip Goerss-Hopkins/obstruct}}
\href{http://hopf.math.purdue.edu/cgi-bin/generate?/Goerss-Hopkins/obstruct}
{\unhbox0}

\bibitem{17a}
\textbf{A Grothendieck}, \emph{Groupes de monodromie en g\'eom\'etrie
  alg\'ebrique I (SGA 7 I)}, Lecture Notes in Mathematics 288, Springer, Berlin
  (1972) \xox{MR}{0354656}

\bibitem{17b}
\textbf{A Grothendieck}, \emph{Groupes de monodromie en g\'eom\'etrie
  alg\'ebrique II (SGA 7 II)}, Lecture Notes in Mathematics 340, Springer,
  Berlin (1973) \xox{MR}{0354657}

\bibitem{16}
\textbf{A Grothendieck}, \emph{Rev\^etements \'etales et groupe fondamental
  ({SGA} 1)}, Documents Math\'ematiques (Paris) 3, Soci\'et\'e Math\'ematique
  de France, Paris (2003) \xox{MR}{2017446}

\bibitem{18}
\textbf{R Hain}, \emph{Lectures on the Hodge--de Rham theory of the fundamental
  group of $\mathbb{P}_{1}-{\{0,1\infty\}}$}
\ Available at \setbox0\hbox{\makeatletter\@url
{http://www.math.duke.edu/~hain/aws/}}
\href{http://www.math.duke.edu/~hain/aws/}
{\unhbox0}

\bibitem{19}
\textbf{M Hazewinkel}, \href{http://dx.doi.org/10.1006/aima.2001.2017}
  {\emph{The algebra of quasi-symmetric functions is free over the integers}},
  Adv. Math. 164 (2001) 283--300 \xox{MR}{1878283}

\bibitem{20}
\textbf{M Hazewinkel}, \href{http://dx.doi.org/10.1007/s10440-004-5635-z}
  {\emph{Symmetric functions, noncommutative symmetric functions and
  quasisymmetric functions. {II}}}, Acta Appl. Math. 85 (2005) 319--340
  \xox{MR}{2128927}

\bibitem{23}
\textbf{M\,E Hoffman}, \href{http://dx.doi.org/10.1006/jabr.1997.7127}
  {\emph{The algebra of multiple harmonic series}}, J. Algebra 194 (1997)
  477--495 \xox{MR}{1467164}

\bibitem{21}
\textbf{M\,J Hopkins}, \emph{Global methods in homotopy theory}, from:
  ``Homotopy theory (Durham, 1985)'', London Math. Soc. Lecture Note Ser. 117,
  Cambridge Univ. Press, Cambridge (1987)  73--96 \xox{MR}{932260}

\bibitem{22}
\textbf{M\,J Hopkins}, \textbf{B Gross}, \emph{The rigid analytic period
  mapping, Lubin--Tate space, and stable homotopy theory}
  \xox{arXiv}{math.AT/9401220}

\bibitem{24}
\textbf{M Hovey}, \textbf{J\,H Palmieri}, \emph{The structure of the
  {B}ousfield lattice}, from: ``Homotopy invariant algebraic structures
  (Baltimore, MD, 1998)'', Contemp. Math. 239, Amer. Math. Soc., Providence, RI
  (1999)  175--196 \xox{MR}{1718080}

\bibitem{25}
\textbf{M Hovey}, \textbf{J\,H Palmieri}, \textbf{N\,P Strickland},
  \emph{Axiomatic stable homotopy theory}, Mem. Amer. Math. Soc. 128 (1997)
  x+114 \xox{MR}{1388895}

\bibitem{26}
\textbf{M Hovey}, \textbf{N\,P Strickland}, \emph{Morava $K$--theories and
  localisation}, Mem. Amer. Math. Soc. 139 (1999) viii+100 \xox{MR}{1601906}

\bibitem{27}
\textbf{M\,M Kapranov}, \emph{Analogies between the {L}anglands correspondence
  and topological quantum field theory}, from: ``Functional analysis on the eve
  of the 21st century, Vol.\ 1 (New Brunswick, NJ, 1993)'', Progr. Math. 131,
  Birkh\"auser, Boston (1995)  119--151 \xox{MR}{1373001}

\bibitem{28}
\textbf{M Kontsevich}, \emph{Operads and motives in deformation quantization}
  \xox{arXiv}{math.QA/990405}

\bibitem{29}
\textbf{H Krause}, \emph{Cohomological quotients and smashing localizations}
  \xox{arXiv}{math.RA/0308291}

\bibitem{30}
\textbf{S Lichtenbaum}, \href{http://dx.doi.org/10.1112/S0010437X04001150}
  {\emph{The {W}eil-\'etale topology on schemes over finite fields}}, Compos.
  Math. 141 (2005) 689--702 \xox{MR}{2135283}

\bibitem{34}
\textbf{J Morava}, \emph{The motivic Thom isomorphism}
  \xox{arXiv}{math.AT/0306151}

\bibitem{32}
\textbf{J Morava}, \emph{The {W}eil group as automorphisms of the
  {L}ubin--{T}ate group}, from: ``Journ\'ees de G\'eom\'etrie Alg\'ebrique de
  Rennes (Rennes, 1978), Vol. I'', Ast\'erisque 63, Soc. Math. France, Paris
  (1979)  169--177 \xox{MR}{563465}

\bibitem{33}
\textbf{J Morava}, \emph{Some {W}eil group representations motivated by
  algebraic topology}, from: ``Elliptic curves and modular forms in algebraic
  topology (Princeton, NJ, 1986)'', Lecture Notes in Math. 1326, Springer,
  Berlin (1988)  94--106 \xox{MR}{970283}

\bibitem{36}
\textbf{A Neeman}, \href{http://dx.doi.org/10.1016/0040-9383(92)90047-L}
  {\emph{The chromatic tower for {$D(R)$}}}, Topology 31 (1992) 519--532
  \xox{MR}{1174255}

\bibitem{37}
\textbf{J Rognes}, \emph{Galois extensions of structured ring spectra}
  \xox{arXiv}{math.AT/0502183}

\bibitem{39}
\textbf{J-P Serre}, \emph{Local class field theory}, from: ``Algebraic Number
  Theory (Proc. Instructional Conf., Brighton, 1965)'', Thompson, Washington,
  D.C. (1967)  128--161 \xox{MR}{0220701}

\bibitem{38}
\textbf{J-P Serre}, \emph{Cohomologie galoisienne}, Lecture Notes in
  Mathematics 5, Springer, Berlin (1994) \xox{MR}{1324577}

\bibitem{40}
\textbf{D Sullivan}, \emph{Geometric topology. {P}art {I}}, Massachusetts
  Institute of Technology, Cambridge, Mass. (1971) \xox{MR}{0494074}

\bibitem{41}
\textbf{J Tate}, \emph{Number theoretic background}, from: ``Automorphic forms,
  representations and $L$--functions (Proc. Sympos. Pure Math., Oregon State
  Univ., Corvallis, Ore., 1977), Part 2'', Proc. Sympos. Pure Math., XXXIII,
  Amer. Math. Soc., Providence, R.I. (1979)  3--26 \xox{MR}{546607}

\bibitem{42}
\textbf{B Toen}, \emph{Vers une interpr\'etation galoisienne de la th\'eorie de
  l'homotopie}, Cah. Topol. G\'eom. Diff\'er. Cat\'eg. 43 (2002) 257--312
  \xox{MR}{1949660}

\bibitem{43}
\textbf{T Torii},
  \href{http://muse.jhu.edu/journals/american_journal_of_mathematics/v125/125.%
5torii.pdf} {\emph{On degeneration of one-dimensional formal group laws and
  applications to stable homotopy theory}}, Amer. J. Math. 125 (2003)
  1037--1077 \xox{MR}{2004428}

\bibitem{44}
\textbf{V Voevodsky}, \emph{Triangulated categories of motives over a field},
  from: ``Cycles, transfers, and motivic homology theories'', Ann. of Math.
  Stud. 143, Princeton Univ. Press, Princeton, NJ (2000)  188--238
  \xox{MR}{1764202}

\bibitem{45}
\textbf{A Weil}, \emph{Basic number theory}, Springer, New York (1974)
  \xox{MR}{0427267}

\end{thebibliography}

\end{document}